\documentclass[a4paper]{amsart}
\usepackage{amsthm,amsmath,amssymb,mathrsfs,url,amsaddr} 
\theoremstyle{definition}
\newtheorem{thm}{Theorem}[section]
\newtheorem*{thm*}{Theorem}
\newtheorem{lem}[thm]{Lemma}
\newtheorem{prop}[thm]{Proposition}
\newtheorem{cor}[thm]{Corollary}

\newtheorem{remark}[thm]{Remark}

\newcommand\DN{\newcommand}
\DN\ts{\times}
\DN\lref[1]{Lemma~\ref{#1}}
\DN\tref[1]{Theorem~\ref{#1}}
\DN\pref[1]{Proposition~\ref{#1}}
\DN\sref[1]{Section~\ref{#1}}
\DN\corref[1]{Corollary~\ref{#1}}
\numberwithin{equation}{section}

\DN\R{\mathbb{R}}
\DN\C{\mathbb{C}}

\DN\xx{\mathbf{x}}
\DN\yy{\mathbf{y}}
\DN\zz{\mathbf{z}}
\DN\ww{\mathbf{w}}
\DN\pp{\mathbf{p}}

\DN\WN{W_{}^N}
\DN\WNp{W_{}^{N+1} }
\DN\intnWN{\mathring{W}_{\ge}^N}
\DN\intWNp{\mathring{W}_{}^{N+1}}
\DN\intnWNp{\mathring{W}_{\ge}^{N+1}}
\DN\WNNp{W^{N,N+1}}
\DN\nWNN{W_{\ge}^{N,N}}
\DN\nWN{W_{\ge}^{N}}
\DN\nWNp{W_{\ge}^{N+1}}

\DN\LNNp{\Lambda_N^{N+1}}
\DN\LaNN{\Lambda_{\alpha, N}^{N}}
\DN\LaNNp{\Lambda_{\alpha, N}^{N+1}}

\begin{document}
\title{The intertwining property for Laguerre processes with a fixed parameter}
\author{ALEXANDER I. BUFETOV }
\address{Steklov Mathematical Institute of RAS, Moscow, Russia, Institute for Information Transmission Problems, Moscow, Russia, CNRS Institut de Math\'ematiques de Marseille, France}
\email{ bufetov@mi.ras.ru}

\author{YOSUKE KAWAMOTO}
\address{Graduate school of environmental, life, natural science and technology, Okayama University, Okayama, Japan}
\email{y-kawamoto@okayama-u.ac.jp}


\begin{abstract}
	We investigate the intertwining of Laguerre processes of parameter $\alpha$ in different dimensions. 
	We introduce a Feller kernel that depends on $\alpha $ and intertwines the  $\alpha$-Laguerre process in $N+1$ dimensions and that in $N$ dimensions.
	When $\alpha $ is a non-negative integer, the new kernel is interpreted in terms of the conditional distribution of the squared singular values: if the singular values of a unitarily invariant random matrix of order $(N+\alpha +1) \times (N+1)$ are fixed, then the those of its $(N+\alpha) \times N $ truncation matrix are given by the new kernel.
\end{abstract}
\maketitle

\section{Introduction and main results}\label{s:1}
\subsection{Intertwining relations for the Laguerre processes}

This note is devoted to intertwining relations for stochastic processes arising in random matrix theory.
More precisely, we set the Weyl chamber $W^N =\{ \mathbf{x} =(x_1,\ldots,x_N)\in \R^N \,;\,  x_1\le \ldots \le x_N  \}$, and for $\xx =(x_1,\ldots, x_{N+1}) \in \WNp $, we introduce the set
\begin{align*}
\WNNp (\xx )& =\{\yy =(y_1,\ldots,y_N) \in W^N \,;\, x_1 \le y_1\le x_2\le \ldots \le y_N \le x_{N+1} \} 
.\end{align*}
For $\yy=(y_1, \ldots,y_N)\in\WN $, let $$\Delta_{N}(\yy)= \prod_{1\le i< j \le N}^{}(y_j-y_i)$$ be the Vandermonde determinant.
For $\xx \in \intWNp := \{ \mathbf{x} \in \WNp \,;\,  x_1 < \ldots < x_{N+1} \}$, let $\LNNp (\xx, \cdot) $ be a probability measure on $W^N$ given by
\begin{align}\label{:21a}
\LNNp (\xx, d\yy ) = N! \cdot \frac{\Delta _N(\yy )}{\Delta _{N+1} (\xx _{})} \mathbf{1}_{ \WNNp (\xx )}(\yy) d\yy 
.\end{align}
Then, $\LNNp $ defines a Feller kernel $W^{N+1} \dashrightarrow W^N $ \cite[Lemma 2.5]{Ass20}. 

Let $T_t^N$ be the Markov semigroup associated with a diffusion on $W^N$. 
We say $T_t^N$ and $T_t^{N+1}$ are intertwined by $\LNNp $ if the relation
\begin{align*}
T_{t}^{N+1} \LNNp =\LNNp T_{t}^{N}
\end{align*}
holds for any $t\ge 0$.
Typical diffusions intertwined by $\LNNp $ include the non-intersecting Brownian motions, also known as the Dyson Brownian motion for $\beta =2$, and its Ornstein-Uhlenbeck counterpart \cite{War07} (see \cite{GoS15,RaS18}  for general inverse temperature $\beta$).
The diffusions introduced in \cite{Ass20}, that leave the Hua-Pickrell measures invariant, are also intertwined by $\LNNp $ .

In this paper, we focus on the intertwining relations concerning Laguerre processes.
For a real number $\alpha >-1$, we consider a diffusion on $[0,\infty)$  associated with the infinitesimal operator
\begin{align}\label{:11f}
L_{\alpha ,x}^{} :=L_\alpha^{} :=x \frac{d ^2}{d x^2}+(-x+1+\alpha ) \frac{d}{dx}
\end{align}
with the following boundary conditions (see, for example, \cite{EtK86} for a detailed discussion of boundary conditions): 
the point $\infty$ is a natural boundary, the origin is an entrance boundary for $\alpha\ge0$ and a regular boundary for $-1<\alpha<0$, in which case we impose the reflecting boundary condition. 

Let $p_{\alpha,t}^{}(x,y):=p_{\alpha}^{}(t,x,y)$ be the transition density of the diffusion associated with $L_\alpha$.
Note that $p_{\alpha,t}$ is the transition function of the solution to the stochastic differential equation
\begin{align*}
dX_t= \sqrt{2X_t} dB_t +(-X_t+\alpha +1) dt
.\end{align*}

We define the non-negative Weyl chamber as $\nWN = W^N \cap[0,\infty)^N $ and let $\intnWN := \{ \mathbf{x} =(x_1,\ldots,x_{N})\in\nWN  \,;\, 0< x_1 < \ldots < x_{N} \}$.
One can direclty check that the Vandermonde determinant $\Delta _N(\mathbf{x})$ is an eigenfunction of the second-order operator $\sum_{1\le i\le N} L_{\alpha ,x_i} $ with eigenvalue $\lambda ^{N}=- N(N-1)/2 $ \cite{Kon08}.
For $(t,\xx, \yy)\in (0,\infty)\times \intnWN \times \nWN $,  consider the Karlin-McGregor transition density of $N$ particles $L_\alpha^{}$-diffusions h-transformed by $\Delta_N(\mathbf{x})$: 
\begin{align*}
\pp_{\alpha ,t}^N(\mathbf{x},\mathbf{y}):=\pp_{\alpha }^N(t, \mathbf{x},\mathbf{y}):=\exp (-\lambda^{N} t) \frac{\Delta_N(\yy ) }{\Delta_N(\xx )} \det_{i,j=1}^{N} [p_{\alpha,t}^{}(x_i,y_j) ].
\end{align*}
Then, $\pp _\alpha ^N$ is the transition density of the non-colliding system of $N$ particles $L_\alpha$-diffusions, which we call the Laguerre process.
By definition, analysis of the Laguerre processes comes down to that of $p_{\alpha, t}$.
To compute the transition density of $L_\alpha$, we will use the following two h-transformations:
\begin{align}\label{:33a}
(\hat{m}_{\alpha}^{})^{} \circ \hat{L} _{\alpha}^{} \circ (\hat{m}_{\alpha}^{})^{-1} &=-1+  L_{\alpha+1}^{} 
\\ \label{:33c}
x^{-\alpha}\circ L_{-\alpha} \circ x^{\alpha} &=-\alpha  +L_{\alpha}
.\end{align}
Here, $\hat{L} _{\alpha}^{}=x (d ^2 /d x^2)+ ( x-\alpha) ( d/ dx)$ is the Siegmund dual of $L_\alpha$ and  $\hat{m}_{\alpha}^{}  (x)=  e^{x} x^{-\alpha -1}$ is its speed measure (see \sref{s:3}).
Both of the two equations are key to deriving the main result \tref{t:21}.

The non-colliding Laguerre process is also obtained as a solution to the stochastic differential equation.
Actually, by the same computation for h-transformation as in \cite{KoO01}, the Laguerre process satisfies
\begin{align}\label{:PND}
dX_t^{i} &=\sqrt{2 X_t^{i} } dB_t^i +\Big( -X_t^{i} +\alpha +1 +\sum_{j\neq i}^N\frac{2 X_t^{i} }{X_t^{i} -X_t^{j} } \Big)dt, \qquad 1\le i\le N 
,\end{align}
where $\{B^i\}_{i=1}^{N}$ is the $N$-dimensional Brownian motion.
This stochastic differential equation has a unique strong solution for any $\xx\in\nWN$ \cite[Thorem 2.2]{GrM14} (we will describe this fact in \lref{l:34} for readers convenience).

The significance of the Laguerre process arises from the random matrix theory.
Specifically, let $m_\alpha^N $ be a probability measure on $\nWN $ given by
\begin{align}\label{:11a}
m_\alpha^N(d \xx ) =\frac{1}{Z_{N,\alpha}} \prod_{1\le i<j\le N}(x_i-x_j)^2 \prod_{k=1}^N x_k^\alpha e^{-x_k} d \xx
,\end{align}
which is referred to as the Laguerre ensemble.
The ensemble \eqref{:11a} gives the distribution of the squared singular values of $(N+\alpha) \times N$ complex Wishart matrix, also known as the Laguerre unitary ensembles. 
The Laguerre process of parameter $\alpha$ leave $m_\alpha^N$ invariant.
Additionally, it is worth noting that the stochastic differential equation \eqref{:PND} is an Ornstein-Uhrenbeck counterpart of the non-colliding squared Bessel processes, which is deeply investigated in the context of random matrix theory \cite{KaT11}.
Actually, the non-colliding squared Bessel process describes the squared singular values process of a rectangular matrices whose entries are independent complex Brownian motions \cite{Bru91,KoO01}.

Let $\{ T_{\alpha,t}^N \}_{t\ge 0}$ be the Markovian semigroup associated with $\pp_{\alpha}^N$.
In \lref{l:34}, we establish the shifted intertwining relation for $\LNNp $ in the sense that  
\begin{align}\label{:34a}
T_{\alpha,t}^{N+1} \LNNp  =\LNNp T_{\alpha +1,t}^{N}
.\end{align}
Shifted intertwining here means that the parameters of Laguerre processes in the left and right hand sides are different.
Remark that a similar shifted intertwining relation is known \cite[Section 3.7]{AOW19}: the relation \eqref{:34a} is an Ornsten-Uhlenbeck counterpart of their result (see also \cite{Ass19} for general $\beta$).
The purpose of this paper is to introduce a new Feller kernel depending on $\alpha$ that intertwines $T_{\alpha,t}^{N+1}$ and $T_{\alpha,t}^{N}$ with the same parameter $\alpha$.

We use the Pochhammer symbol $(x)_{n} =x(x+1)\cdots(x+n-1)$.
Suppose $\alpha >-1 $.
For $\xx\in \mathring{W}_{\ge}^{N+1}$, define $\LaNNp (\xx, \cdot )$ as a measure on $\nWN$ given by
\begin{align}\label{:21b}
\LaNNp (\xx, d\yy) =N! (\alpha+1)_N  \frac{\Delta _N ( \yy )}{\Delta _{N+1} (\xx _{})} \prod_{k=1}^{N} \bigg( \mathbf{1}_{[x_{k-1} , x_{k+1}]}( y_k)  \int_{x_k \vee y_k }^{x_{k+1} \wedge y_{k+1} } \frac{y_k^{\alpha} }{z ^{\alpha+1}} d z  \bigg) d\yy 
.\end{align}
Here, we use the symbol $x_0 =0 $ and $y_{N+1}=\infty$ for notational convenience.
While the definition of \eqref{:21b} is valid for $\xx\in\intnWNp $, we will see in \sref{s:22} that $\LaNNp $ is extended to a Feller kernel $\nWNp \dashrightarrow \nWN $.
It is important to emphasise that \eqref{:21b} is formulated for real parameter $\alpha$.
In the particular case when $\alpha$ is a non-negative integer, the kernel $\LaNNp $ can be interpreted in terms of the radial parts of random matrices.
Actually, we shall see in \tref{t:12} that the kernel $\LaNNp (\xx, \cdot)$ provides the radial parts distribution of a truncation of unitary invariant matrix whose radial parts correspond with $\xx \in \nWNp $.

Setting $C_{\infty} (\nWN )$ to be the set of all continuous functions on $\nWN $ vanishing at infinity, now we state our first main result.
\begin{thm}\label{t:21}
Assume that  $\alpha >-1$.
Then, for any $N\in\mathbb{N} $, $f\in C_{\infty} (\nWN )$, and $t \ge 0$, we have the intertwining relation
\begin{align}\label{:11b}
T_{\alpha,t}^{N+1} \LaNNp f=\LaNNp T_{\alpha,t}^{N}f
.\end{align}
\end{thm}
In contrast to \eqref{:34a}, the parameters of the left and right hand sides in \eqref{:11b} are the same.
Naturally, we find that a sequence of $\alpha$-Laguerre processes $\{T_{\alpha}^N \}_{N\in\mathbb{N} } $ is coherent with kernels $\{\LaNNp \}_{N\in\mathbb{N} }$.

Intertwining relations allow one to construct a limit process on the boundary of branching graphs corresponding to a coherent family of finite-dimensional processes.
Borodin and Olshanski \cite{BoO12} first introduced the method of intertwiners, constructing a Feller process on the boundary of the Gelfand-Tsetlin graph, which describes the branching of irreducible representations of the chain of unitary groups.
Borodin and Olshanski \cite{BoO13, Ols16} and Cuenca \cite{Cue18} applied this approach for branching graphs related to other groups.
For the background on the intertwining property, see ``The Gelfand-Tsetlin graph and Markov processes'' by Olshanski (arXiv:1404.3646).
In their frameworks, the state spaces of finite-dimensional processes are discrete.
Assiotis \cite{Ass20} applied the method of intertwiners in cases where processes in finite dimensions are continuous.
The coherent family coming from the intertwining relation \eqref{:11b} also yields a limit process on a boundary as $N \to \infty $, and this will be done in the sequel to this note.

\subsection{Interpretation of Markov kernels in terms of radial parts of random matrices}

Hereafter in this subsection, $\alpha$ is supposed to be a non-negative integer.
In this case, the Markov kernel $\LaNNp $ has the following interpretation in the context of random matrix theory.

Let $M_{m,n}(\C)$ be the space of $m\times n$ matrices with complex entries, and  for brevity write $M_{n}(\C)=M_{n,n}(\C)$.
Introduce the following subsets  $H_{n}(\C), \mathbb{U}(n) \subset M_{n}(\C) $:  $H_{n}(\C) $ is the space of Hermite matrices of order $n$, and $\mathbb{U}(n)$ is the space of unitary matrices of order $n$.

For $m_1\ge m_2, n_1\ge n_2 $, let $\pi_{m_2, n_2}^{m_1, n_1} : M_{m_1, n_1}  (\C) \to M_{m_2,n_2} (\C) $ be the natural projection sending an $m_1 \times n_1 $ matrix to its upper left $m_2 \times n_2$ corner.
We employ the expression $\pi_{m_2,n_2}^{m_1}$ in place of $\pi_{m_2,n_2}^{m_1,n_1}$ if $m_1=n_1$, and use a similar symbol for $m_2=n_2$.

We define a map $\mathfrak{eval}_{n} : H_{n}(\C) \to W^n  $ as 
\begin{align*}
\mathfrak{eval}_{n} (X)=(\lambda_{1}(X), \ldots , \lambda_{n}(X))
,\end{align*}
where $(\lambda_{i}(X))_{i=1}^n$ is the eigenvalues of $X$ arranged in non-decreasing order.
Furthermore, define the radial part $\mathfrak{rad}_{n} : M_{m,n} (\C) \to W_{\ge}^n $ as $\mathfrak{rad}_{n} (X)=\mathfrak{eval} _n(X^{*}X)$.
Let a probability measure $P_{\mathfrak{eval}}^n [X ] $ on $W_{\ge}^n $ be the distribution of the eigenvalues of a random matrix $X \in M_{n}(\C) $.
Similarly, let a probability measure $P_{\mathfrak{rad}}^n [X]  $ on $W_{\ge}^n$ denote the distribution of the radial part of a random matrix $X\in M_{m,n}(\C) $.

Let $U_{N+1}\in \mathbb{U} (N+1)$ be a Haar distributed random matrix and let $\mathrm{diag}(x_1,\ldots , x_{N+1} ) $ be the square matrix of order $N+1$ with deterministic diagonal elements $x_1,\ldots, x_{N+1}$.
Then, it is well known that $\LNNp (\xx , \cdot) $ is given by the formula
\begin{align*}
\LNNp (\xx , \cdot)=P_{\mathfrak{eval}}^{N} \big[\pi_{N}^{N+1} ( U_{N+1}^{*} \mathrm{diag}(x_1,\ldots , x_{N+1} ) U_{N+1} ) \big]
\end{align*}
valid for any $\xx\in\WNp$ \cite[Proposition 4.2]{Bar01}. 
Hence, if a random matrix $X_{N+1}\in H_{N+1}(\C)$ is $\mathbb{U}(N+1) $-invariant by conjugation in the sense that $U_{N+1}^{*} X_{N+1} U_{N+1} \stackrel{\mathrm{law}}{=} X_{N+1}$ for any unitary matrix $U_{N+1} \in \mathbb{U}( N+1) $, then we have the relation
\begin{align}\label{:31c}
P_{\mathfrak{eval}}^{N+1}[X_{N+1}] \LNNp=P_{\mathfrak{eval}}^{N}[ \pi_{N}^{N+1}( X_{N+1}) ]
.\end{align}

The next theorem gives a similar expression for the Markov kernel $\LaNNp $, and consequencty, a relation between the radial parts of a unitary invariant rectangular matrix and  its truncation.
A random matrix $X_{m,n}\in M_{m,n}(\C)$ is said to be $\mathbb{U}(m) \times \mathbb{U}(n) $-invariant if $X_{m,n}\stackrel{\mathrm{law}}{=} V_{m} X_{m,n} U_{n}$ for any fixed unitary matrices $V_{m}\in \mathbb{U}(m), U_{n}\in \mathbb{U}(n)$.
\begin{thm}\label{t:12}
Let $\alpha$ be a non-negative integer.
For a $\mathbb{U}(N+\alpha+1) \times \mathbb{U}( N+1)$-invariant random matrix $X_{N+\alpha+1,N+1} \in M_{N+\alpha+1,N+1}(\C)$, let $X_{N+\alpha, N}=\pi_{N+\alpha,N}^{N+\alpha+1, N+1}( X_{N+\alpha+1,N+1})$ be its truncation.
We then have
\begin{align*}
P_{\mathfrak{rad}}^{N+1} [ X_{N+\alpha+1,N+1} ] \LaNNp =P_{\mathfrak{rad}}^{N} [X_{N+\alpha, N}]
.\end{align*}
\end{thm}

\begin{cor}\label{c:16}
Let $V_{N+\alpha+1}  \in \mathbb{U}(N+\alpha+1) $ and $U_{N+1}  \in \mathbb{U}(N+1) $ be Haar distributed independent random matrices and $D_{N+\alpha+1, N+1} \in M_{N+\alpha+1, N+1}(\C)$ be a deterministic matrix given by
\begin{align*}
D_{N+\alpha+1,N+1}=
\begin{bmatrix}
	\mathrm{diag}(\sqrt{x_1},\ldots,\sqrt{x_{N+1}} )  \\
\mathbf{0}_{\alpha \times (N+1)}
\end{bmatrix}
\text{ for }\xx =(x_1,\ldots, x_{N+1})\in \nWNp
.\end{align*} 
Then, the Markov kernel $\LaNNp(\xx , d \yy) $ is the same as 
\begin{align*}
P_{\mathfrak{rad}}^{N} \big[ \pi _{N+ \alpha ,N}^{N+\alpha+1, N+1}(V_{N+\alpha+1}  D_{N+\alpha+1, N+1}  U_{N+1}) \big]
.\end{align*}
\begin{proof}
Since $V_{N+\alpha+1}  D_{N+\alpha+1, N+1}  U_{N+1}$ is a $\mathbb{U}(N+\alpha+1) \times \mathbb{U}( N+1)$-invariant random matrix and $P_{\mathfrak{rad}}^{N+1} [V_{N+\alpha+1}  D_{N+\alpha+1, N+1}  U_{N+1}] =\delta_{\xx} $ holds, \tref{t:12} implies the statement of this corollary.
\end{proof}
\end{cor}

The Laguerre ensemble $m_\alpha^{N} $ defined in \eqref{:11a} is identical to the distribution of radial parts of the $(N+\alpha)\times N $ complex Wishart matrix (see, for example, \cite{For10}).
Hence, the relation 
\begin{align}\label{:26b}
m_\alpha^{N+1} \LaNNp =m_\alpha^{N}
\end{align}
immediately follows from \tref{t:12} for non-negative integer $\alpha $.
The equation \eqref{:26b} for general $\alpha $ will be shown in \sref{s:4} from the result of the intertwining relation for the Laguerre processes.

\subsection{Organisation of this paper}
The present paper is organised as follows.
\sref{s:2} is devoted to examine the new kernel $\LaNNp $.
In particular, we present useful expressions of $\LaNNp $, and show its Feller property.
In \sref{s:3}, we will show two shifted intertwining relations, which lead to \tref{t:21}.
We give the proof of \tref{t:12} in \sref{s:4}.

\section{Feller kernels}\label{s:2}
\subsection{Kernels and semigroups}
We first recall concepts of kernels and semigroups on the Euclidean spaces (see, for instance, \cite[section 2]{BoO12} for a more comprehensive discussion).
Let $E$ and $E'$ be Borel subsets of Euclidean spaces.
Consider a function $\Lambda(x,A)$, where $x\in E $ and $A$ is a Borel subset of $E'$.
We say that $\Lambda $ is a Markov kernel from $E$ to $E'$, denoted as $\Lambda : E \dashrightarrow E' $, if the following two conditions hold:
\begin{itemize}
\item
 $\Lambda (x ,\cdot )$ is a Borel probability measure on $E'$ for any $x\in E$.

\item
$\Lambda(\cdot, A) $ is a Borel function on $E$ for any Borel subset $A\in E'$.
\end{itemize}

Let $\mathcal{B}(E) $ and $\mathcal{B}(E') $ denote the Banach spaces of $\mathbb{R} $-valued bounded measurable functions with the sup-norm on $E$ and $E'$, respectively.
Then, a Markov kernel $\Lambda : E \dashrightarrow E' $ defines a linear operator $\mathcal{B}(E') \to \mathcal{B}(E)$ as $(\Lambda f)(x)=\int_{E'}\Lambda (x, dy) f(y)$.
Let $C_\infty(E)$ be the set of all continuous functions on $E$ vanishing at infinity.
We say a Markov kernel $\Lambda : E \dashrightarrow E' $ is Feller if the induced map $\Lambda :\mathcal{B}(E') \to \mathcal{B}(E)$ satisfies $\Lambda(C_\infty(E')) \subset C_\infty(E) $.

The symbols $\mathcal{M}_p(E)$ and $\mathcal{M}_p(E')$ are used to represent the set of all of probability measures on $E$ and $E'$, respectively.
For $\mu \in \mathcal{M}_p(E)$, we set $\mu \Lambda (\cdot)=\int_{E}\Lambda (x,\cdot ) \mu (dx) \in \mathcal{M}_p(E')$.
Thus, a Markov kernel $\Lambda : E \dashrightarrow E' $ induces a linear map $\mathcal{M}_p(E) \to \mathcal{M}_p(E')$.

A Markov semigroup on $E$ is, by definition, a family of Markov kernels $\{ T_t\}_{t \ge 0}$ on $E$ such that $T_0 =1$ and $T_t T_s=T_{t+s}$.
This semigroup induces a semigroup on $\mathcal{B}(E)$.
A Markov semigroup on $E$ is said to be Feller-Dynkin if the following two conditions hold:

\begin{itemize}
\item 
The induced map $T_t: \mathcal{B}(E) \to \mathcal{B}(E)$ satisfies $T_t ( C_\infty(E) ) \subset C_\infty(E) $ for any $t\ge 0$.

\item
The function $t\to T_t$ is strongly continuous in the sense that relation $\lim_{t\to 0} T_t f =f $ in $C_\infty(E)$ holds for any $f\in C_\infty(E)$. 
\end{itemize}

We note that Feller-Dynkin semigroups are often called simply Feller.
To avoid confusion, we distinguish between these two in this paper.
For us, the Feller property means just the property of mapping continuous functions vanishing at infinity into itself; Feller-Dynkin is used if, in addition to this, the strong continuity holds.

\subsection{Feller property of the Markov kernels $\LaNNp $}\label{s:22}

In this section, by the Feller property of Markov kernels we simply mean that the Markov operator maps continuous functions vanishing at infinity into itself.
The Feller property of $\LaNNp $ is demonstrated in \pref{p:24}, which shall be used to establish \tref{t:21}.

The explicit expression of $\LaNNp $ given by \eqref{:21b} is not convenient for computation.
We decompose $\LaNNp $ into two simpler kernels, $\LNNp $ given in \eqref{:21a} and $\LaNN $ introduced below.
This decomposition proves to be useful for establishing the Feller property and proving \tref{t:21}.

For $\zz =(z_1,\ldots, z_{N}) \in \nWN$, we set
\begin{align*}
\nWNN (\zz )& =\{\yy =(y_1,\ldots,y_N) \in \nWN \,;\, 0\le  y_1 \le z_1\le y_2\le \ldots \le y_N \le z_{N} \} 
.\end{align*}
For $\zz \in \intnWN $, we introduce a kernel  $\LaNN (\zz, d\yy) $  on $\nWN $ by setting 
\begin{align}\label{:32c}
\LaNN (\zz, d\yy) = (\alpha+1)_{N} \bigg( \prod_{k=1}^{N}\frac{ y_k^\alpha }{ z_k^{\alpha+1}  }\bigg) \frac{\Delta_N (\yy)}{\Delta_N (\zz )}  \mathbf{1}_{\nWNN (\zz )}(\yy) d\yy
.\end{align}

Hereafter in this subsection, and notably in the proof of \lref{l:31} and \lref{l:32}, we set $z_0=0$ and $z_0^n=0$ for notational convenience. 

\begin{lem}\label{l:31}
If $\alpha >-1$ and $\zz \in \intnWN $, then $\LaNN (\zz, \cdot) $ is a probability measure on $\nWN $.
\begin{proof}
The following computation yields the statement of this lemma:
\begin{align*}
\int_{\nWNN (\zz )} \Big( \prod_{k=1}^{N} y_k^\alpha  \Big)\Delta_{N}(\yy)  d\yy&=\int_{z_0}^{z_1} dy_1 \int_{z_1}^{z_2} dy_2 \ldots \int_{z_{N-1}}^{z_N} dy_N \det_{i,j=1}^{N}[y_i^{j-1+\alpha}]
\\&
= \det_{i,j=1}^{ N} \Big[ \int_{z_{i-1}}^{z_i } y^{j-1+\alpha}dy \Big]
\\&
=\frac{\big(\prod_{k=1}^{N} z_k^{\alpha+1} \big) \Delta_{N}(\zz ) }{(\alpha+1)_N}
.\end{align*}
The last equality results from successive addition of row $i$ to $i+1$ starting from $i=1$.
\end{proof}
\end{lem}

The explicit density \eqref{:32c} gives the definition of $\LaNN (\zz, d \yy)$ for $\zz\in\intnWN$.
For  $\zz\in \nWN$, the transition kernel $\LaNN (\zz, d \yy)$ is defined by taking the weak limit verified in the following lemma.
For $\zz \in \intnWN$, we write down $\LaNN f $ explicitly as
\begin{align}\label{:32f}
\LaNN f (\zz ) = (\alpha+1)_N  \int_{z_0}^{z_1}\ldots \int_{z_{N-1}}^{z_{N}}\Big( \prod_{k=1}^{N} \frac{y_k^\alpha}{ z_k^{\alpha+1} }  \Big) \frac{\Delta_{N}(\yy)}{\Delta_{N}(\zz )} f(\yy) d \yy
.\end{align}
Note that $\LaNN f $ can be defined only on $\intnWN $ at this stage.

\begin{lem}\label{l:32}
Let $\alpha >-1 $.
Then, for any bounded continuous function $ f \in C_b(\nWN )$, the function $\LaNN f  $ is continuous on $\intnWN $.
Furthermore, it can be continuously extended to $ \nWN$.
\begin{proof}

The continuity on $\intnWN $ immediately follows from \eqref{:32f} by the bounded convergence theorem since the integrand has no singularities.
We will show that, for any $\zz \in \nWN \setminus \intnWN$ and a sequence $\zz ^{n}=(z_1^{n},\ldots,z_N^{n}) \in \intnWN $ satisfying $\zz:=\lim_{n\to \infty} \zz^n $, there exists a limit of $\LaNN f (\zz^n )$ as $n\to\infty$ and it is independent of a choice of an approximating sequence $\{ \zz^n \}_{n\in\mathbb{N} }$.
We first assume that a limit point $\zz=(z_1,\ldots,z_N) $ satisfies
\begin{align}\label{:32g}
z_0 <\ldots <z_{L-1} < z_{L}=\ldots=z_{M} <z_{M+1}<\ldots<z_{N}
\end{align}
for some $0\le L<M\le N$.
Denote $a =z_{L}=\ldots=z_{M} $ (remark that $L=0$ implies $a=0$ because of our notation).

We set $ \yy_{L}^{M}=(y_1,\ldots, y_{L}, y_{M+1},\ldots ,y_N) \in W_{\ge}^{N-M+L}$ and $ \yy_{L}^{M}(a) =(y_1,\ldots, y_{L}, a,\ldots, a, y_{M+1},\ldots ,y_N)\in \nWN$. 
By a straightforward computation similar to in \lref{l:31}, we obtain
\begin{align*}
\int_{z_{L}^n}^{z_{L+1}^n}\ldots \int_{z_{M-1}^n}^{z_{M}^n} \Big( \prod_{k=L+1}^{M} y_k^\alpha \Big) \Delta_{N}(\yy)  \prod_{k=L+1}^{M} d y_k
\\
 = \Big( \prod_{k=L+1}^{M} (z_k^n)^{\alpha +1}\Big) \Delta_{M-L} (z_{L+1}^n,\ldots, z_{M}^n) H(\yy_L^M, \zz ^n)
,\end{align*}
where $H$ is a continuous function.
Hence, the kernel $G$ given by the formula
\begin{align*}
G(\yy_L^M, \zz ^n ):=\int_{z_{L}^n}^{z_{L+1}^n}\ldots \int_{z_{M-1}^n}^{z_{M}^n} \Big( \prod_{k=1}^{N} \frac{y_k^\alpha}{ (z_k^n)^{\alpha+1} }  \Big) \frac{\Delta_{N}(\yy)}{\Delta_{N}(\zz ^n)} \prod_{k=1}^{M} d y_k  
\end{align*}
has no singularity in $\zz ^n$ and quantifiers on $L,M$, the limit $\lim_{n\to \infty} G(\yy_L^M, \zz^n) $ exists by \eqref{:32g}.
Because of \eqref{:32g} and the continuity of $f$ on $\nWN $, for any $\varepsilon >0$, we have
\begin{align}\label{:32h}
\limsup_{n\to \infty} \sup_{\yy \in W_{\ge}^{N,N}(\zz ^n) } |f(\yy)-f( \yy_{L}^{M}(a))| <\varepsilon
.\end{align}
We write
\begin{align*}
F(\yy_{L}^{M} , \zz^n) = \int_{z_{L}^n}^{z_{L+1}^n}\ldots \int_{z_{M-1}^n}^{z_{M}^n} \Big( \prod_{k=1}^{N} \frac{y_k^\alpha}{ (z_k^n)^{\alpha+1} }  \Big) \frac{\Delta_{N}(\yy)}{\Delta_{N}(\zz ^n)}  f(\yy) \prod_{k=L+1}^{M} d y_k
.\end{align*}
A direct computation using inequality \eqref{:32h} gives
\begin{align*}
|F( \yy_{L}^{M} ,\zz^n ) -G(\yy_{L}^{M} ,\zz^n) f( \yy_{L}^{M}(a) ) | \le G(\yy_L^M, \zz^n) \varepsilon
,\end{align*}
which implies $\lim_{n\to\infty} F( \yy_{L}^{M}, \zz^n ) =G(\yy_{L}^{M} ,\zz) f( \yy_{L}^{M}(a))$.
Furthermore, by the same argument for $F$ with the dominated convergence theorem applied to the integral
\begin{align*}
\LaNN f (\zz ^n) =\int_{z_{0}^n}^{z_{1}^n}\ldots \int_{z_{L-1}^n}^{z_{L}^n} \int_{z_{M}^n}^{z_{M+1}^n} \ldots \int_{z_{N-1}^n}^{z_{N}^n}F(\yy_{L}^{M}, \zz^n  )  \prod_{\substack{1\le k \le L \\ M+1\le k\le N}} d y_k
,\end{align*}
we obtain the existence of $\lim_{n\to\infty} \LaNN f (\zz ^n)$, and the limit is clearly independent of the specific choice of sequence $\{ \zz^n\}_{n\in\mathbb{N} }$.
Thus, we have extended $\LaNN f (\zz )$ continuously to the points satisfying \eqref{:32g}.
The proof of the general case is the same.

\end{proof}
\end{lem}

\begin{lem}\label{l:17}
Suppose $\alpha >-1 $.
Then, for any $\xx \in \intnWNp  $, we have
\begin{align*}
\LaNNp (\xx, d\yy )=\LNNp \LaNN (\xx,d\yy) 
.\end{align*}
\begin{proof}
From \eqref{:21a} and \eqref{:32c}, a direct computation gives
\begin{align*}&
\int \LNNp (\xx, d \zz) \LaNN (\zz, d \yy) 
\\&
=N! (\alpha+1)_N  \frac{\Delta _N ( \yy )}{\Delta _{N+1} (\xx _{})}  d\yy  \int \bigg( \prod_{k=1}^{N} \frac{y_k^{\alpha}}{z_k ^{\alpha+1}} \bigg) \mathbf{1}_{ \WNNp (\xx )}(\zz ) \mathbf{1}_{\nWNN (\zz )}(\yy) d\zz 
\\&
=N! (\alpha+1)_N  \frac{\Delta _N ( \yy )}{\Delta _{N+1} (\xx _{})}  d\yy \prod_{k=1}^{N} \bigg( \int_{x_k \vee y_k }^{x_{k+1} \wedge y_{k+1} }   \frac{y_k^{\alpha} }{z ^{\alpha+1}} d z  \mathbf{1}_{\{ x_{k} \vee y_{k} \le x_{k+1} \wedge y_{k+1} \}}  \bigg)
\\&
=\LaNNp (\xx, d\yy)
\end{align*}
by the definition \eqref{:21b}.
\end{proof}
\end{lem}

Take $\xx \in \nWNp $ and $\xx^n \in \intnWNp $ satisfying $\lim_{n\to\infty }\xx^n =\xx$.
Then, we have 
\begin{align*}
 \lim_{n\to\infty }\LNNp (\xx^n ,\cdot )=\LNNp (\xx ,\cdot)
\end{align*}
weakly \cite[Lemma 2.5]{Ass20}.
Therefore, for any $f\in C_{b}(\nWN)$, we obtain 
\begin{align*}
\lim_{n\to\infty} \LNNp \LaNN f (\xx^n)= \LNNp \LaNN f (\xx )
.\end{align*}
Combining this with \lref{l:17}, the function $\LaNNp f$ can be continuously extended to $\nWN$.
Thus, we can define the kernel $\LaNNp $ on $\nWN$.
Accordingly, the decomposition in \lref{l:17} is extended as follows:

\begin{prop}\label{p:23}
Suppose $\alpha>-1$.
Then, for any $\xx \in \nWNp  $, we have
\begin{align*}
\LaNNp (\xx, d\yy)=\LNNp \LaNN (\xx, d\yy) 
.\end{align*}
\end{prop}

\lref{l:31} and \lref{l:32} imply that $\LaNN $ is a Markov kernel $\nWN \dashrightarrow \nWN$.
Recall that $\LNNp $ is Markov.
Therefore, $\LaNNp $ is also a Markov kernel $\nWNp \dashrightarrow \nWN $ from \pref{p:23}.
We now show that $\LaNN $ and $\LaNNp $ are Feller kernels.

\begin{prop}\label{p:24}
Suppose $\alpha >-1$.
Then, the following (i) and (ii) hold.
\begin{itemize}
\item[(i)]
$\LaNN f \in C_{\infty}(\nWN)$ for any $f\in  C_{\infty}(\nWN )$.

\item[(ii)]
$\LaNNp f \in C_{\infty}(\nWNp )$ for any $f\in  C_{\infty}(\nWN )$.
\end{itemize}
\begin{proof}

The statement (ii) immediately follows from (i) and \pref{p:23} because
\begin{align}\label{:24o}
\LNNp f\in C_\infty (\WNp )
\end{align}
holds for any $f\in C_\infty( \WN )$ \cite[Lemma 2.5]{Ass20}.
We proceed to the proof of (i).
The continuity of $\LaNN f$ on $\nWN$ is proved in \lref{l:32}, and to establish (i) it remains to show the relation
\begin{align}\label{:24n}
\lim_{n\to\infty}\LaNN f(\zz ^n)=0
\end{align}
for any sequence $\{ \zz^n \}_{n \in\mathbb{N} } \subset \nWN $ satisfying $\lim_{n\to\infty } \zz^n =\infty$.

Suppose $\zz \in \intnWN$.
Applying the mean value theorem to integrals for $y_N,\ldots, y_2$ successively, we see that there exists $\ww  =(w_2,\ldots, w_N)\in \mathring{W}_{\ge}^{N-1}$ satisfying 
\begin{align}\label{:24m}
z_1 < w_2 <z_{2}< \ldots < z_{N-1} < w_{N} < z_{N}
\end{align}
such that  
\begin{align*}
\int_{z_1}^{z_2}\!\!\!\!\! \ldots\!\! \int_{z_{N-1}}^{z_{N}} \!\!\!\Big( \prod_{k=1}^{N} y_k^\alpha \Big) \Delta_{N}(\yy)  f(\yy) \! \prod_{k=2}^N d y_k= y_1^\alpha \Big( \prod_{k=2}^{N} w_k ^\alpha\Big) \Delta_{N}(y_1, \ww ) f(y_1, \ww )  \! \prod_{k=2}^{N}(z_k-z_{k-1})  
.\end{align*}
Thereby, we obtain
\begin{align}\label{:24b}
| \LaNN f (\zz ) |
\le  (\alpha+1)_N   \frac{ \Big( \prod \limits_{k=2}^{N} w_k^\alpha \Big)  \prod \limits_{k=2}^{N}(z_k-z_{k-1}) }{ \Big( \prod \limits_{k=1}^{N}  z_k^{\alpha+1} \Big)\Delta_{N}(\zz)} \int_{0}^{z_{1}}  y_{}^\alpha \Delta_{N}(y_{}, \ww ) | f(y_{}, \ww ) |dy_{}  
.\end{align}
From \eqref{:24m}, it holds that 
$$ \Delta_{N}(y_{}, \ww ) \le \Delta_{N-1}(\ww) \prod_{k=2}^N w_k
$$ for $y_{}\le z_1$.
Combining this with \eqref{:24b} yields
\begin{align}\label{:24j}
| \LaNN f (\zz ) |
&
\le (\alpha+1)_N  \frac{ \Big( \prod \limits_{k=2}^{N} w_k^{\alpha+1} \Big) \prod\limits_{k=2}^{N}(z_k-z_{k-1}) \Delta_{N-1}(\ww)  }{ \Big( \prod \limits_{k=1}^{N}  z_k^{\alpha+1} \Big)\Delta_{N}(\zz)} {\displaystyle \int _{0}^{z_{1}}  y_{}^\alpha  | f(y_{}, \ww )| dy } 
\\&\notag
= (\alpha+1)_N \bigg( \prod_{k=2}^{N} \frac{ w_k }{  z_k }\bigg)^{\alpha+1} \bigg( \prod_{\substack{1\le i<j\le N\\  j\neq i+1}} \frac{w_j-w_{i+1}}{z_j-z_i} \bigg)  \frac{\displaystyle{\int_{0}^{z_{1}}  y_{}^\alpha  | f(y_{}, \ww )| dy}}{z_1^{\alpha+1}}
\\& \notag
< (\alpha+1)_N \frac{w_N-w_{N-1}}{z_{N}-z_{N-2}}  \frac{{\displaystyle \int_{0}^{z_{1}}  y_{}^\alpha  | f(y_{}, \ww )| dy}}{z_1^{\alpha+1}}
\end{align}
Here, the last inequality results from \eqref{:24m}.
Actually, $(\prod_{k=2}^{N} w_k/z_k) ^{\alpha+1}  < 1$ holds since $w_k < z_k$  and  $\alpha >-1$, and 
\begin{align}\label{:24g}
\frac{ w_j-w_{i+1}  }{ z_j-z_i } < 1
\end{align}
also holds from interlace relation $z_{i}^{}< w_{i+1}< w_{j} < z_j  $ for $i+1 <j$.

For any $\varepsilon >0$, there exists $p$ such that $ | f(y_{}, \ww )| <\varepsilon $ for $w_N >p $.
Then, if  $w_N >p $, the inequalities \eqref{:24j} and \eqref{:24g} yield
\begin{align*}
| \LaNN f (\zz ) |&< (\alpha+1)_N  \frac{{\displaystyle \int_{0}^{z_{1}} } y_{}^\alpha \varepsilon dy_{}}{z_1^{\alpha+1}} = \varepsilon (\alpha+2)_N 
.\end{align*}
On the other hand, if $w_N \le p $, then \eqref{:24j} with $z_{N-2}<w_N$ implies the inequality
\begin{align*}
| \LaNN f (\zz ) |
&
\le  (\alpha+1)_N \frac{p}{z_{N}-p}  || f|| \frac{{\displaystyle \int_{0}^{z_{1}}  y_{}^\alpha   dy_{} }}{z_1^{\alpha+1} }
=  (\alpha+2)_N || f||  \frac{p}{z_{N}-p}
,\end{align*}
where $||f||$ is the sup-norm on $\nWN$. 
In either case, we then see $| \LaNN f (\zz ) |$ is sufficiently small for sufficiently large $z_N \to \infty$.
Thus, observing that $\lim_{n\to\infty}\zz^n =\infty$ implies $\lim_{n\to\infty }z_N^n=\infty$, we have proved \eqref{:24n} if $\zz^n \in \intnWN $.
The relation \eqref{:24n} for $\zz^n \in \nWN $ follows from the continuity of $\LNNp f $ on $\nWN$ by \lref{l:32}.
The proof of (i), and hence of \pref{p:24} is complete.

\end{proof}
\end{prop}

\section{Intertwining for the Laguerre processes}\label{s:3}
\subsection{The Laguerre processes}
Here we collect some facts about Laguerre processes.
We first recall the following fact introduced in \sref{s:1}, although we will not use this result in this paper. 
\begin{lem}\label{l:35}
The stochastic differential equation \eqref{:PND} has a unique strong solution for any starting point $\xx\in\nWN$.
\begin{proof}
Because the stochastic differential equation \eqref{:PND} becomes 
\begin{align*} 
dX_t^{i} &=\sqrt{2 X_t^{i} } dB_t^i +\Big( -X_t^{i} +\alpha +N +\sum_{j\neq i}^{N} \frac{ X_t^{i}+X_t^{j} }{X_t^{i} -X_t^{j} } \Big)dt,
\end{align*}
by an easy computation, the strong uniqueness results from \cite[Theorem 2.2]{GrM14} by taking $H_{ij}(x,y)=x+y $, $\sigma_{i}(x)=\sqrt{2x}$, and $b_i(x)=-x+\alpha+N$ (see also \cite[Corollary 6.4]{GrM14}). 
\end{proof} 
\end{lem}

We next quote the Feller-Dynkin property of the Laguerre process.
\begin{lem}\cite[Proposition 1.7]{Ass19b} \label{l:33}
The semigroup $T_{\alpha,t}^{N}  $ is Feller-Dynkin in the sense that, for any $f\in C_{\infty } (\nWN )  $, we have $T_{\alpha, t}^{N}  f\in C_{\infty}  (\nWN )$ for any $t\ge 0 $ and $\lim_{t\to 0} T_{\alpha, t}^{N}  f =f $.
\end{lem}

\subsection{Proof of \tref{t:21} }
We consider the operator $L_{\alpha}^{}$ in \eqref{:11f} for any $\alpha\in\R$.
When $\alpha \le -1$, the origin is an exit boundary.
Let $\hat{L}_{\alpha}^{}$ be the Siegmund dual operator of $L_{\alpha}^{}$, that is, 
\begin{align*}
\hat{L}_\alpha^{}=x \frac{d ^2}{d x^2}+ ( x-\alpha) \frac{d}{dx}
.\end{align*}
Here, the origin is an exit boundary for $\alpha\ge0$, a regular absorbing boundary for $-1<\alpha<0$, and an entrance boundary for $\alpha \le -1$.
The speed measure for $\hat{L}_\alpha^{} $ is given by $\hat{m}_\alpha (x)=  e^{x} x^{-\alpha -1}$.
Here, we recall that the scale function $s$ and the speed measure $m$ for $a(x) ( d^2/dx^2) +b(x) (d/dx)$ are given by
\begin{align*}
s'(x)=\exp \bigg(- \int_c^x \frac{b(y)}{a(y)}dy \bigg), \qquad m(x)= \frac{1}{s'(x) a(x)}
\end{align*}

We need two types of Doob's h-transformation.
\begin{lem}
\begin{itemize}
\item[(i)]
For any $\alpha \in\R$, we have
\begin{align}\label{:21c}
e^{ t} \hat{p}_{{\alpha } ,t}^{} (x,y)\frac{(\hat{m}_{\alpha }^{}(y))^{-1}} {(\hat{m}_{\alpha}^{}(x) )^{-1} }= p_{\alpha+1,t} (x,y)
\end{align}

\item[(ii)]
For any $\alpha >-1$, we have
\begin{align} \label{:32a}
e^{ \alpha t} p_{-\alpha,t}(x,y) \frac{y^{\alpha } }{x^{\alpha }} =p_{\alpha,t}(x,y)
.\end{align}
\end{itemize}
\begin{proof}

We see that a positive function $(\hat{m}_\alpha (x) )^{-1} = e^{-x} x^{\alpha +1}$ is  an eigenfunction of $\hat{L}_\alpha $ with eigenvalue $-1$. 
Then, we see \eqref{:33a} from a straightforward calculation, which implies \eqref{:21c}.
Furthermore, note that a positive function $x^{\alpha } $ is an eigenfunction of $L_{-\alpha}$ with eigenvalue $- \alpha $. 
A direct computation yields \eqref{:33c}, which implies \eqref{:32a}
\end{proof}
\end{lem}

We define positive kernels $\hat{\Lambda}_{\alpha, N}^{N+1} : \nWNp \dashrightarrow \nWN $ and $\hat{\Lambda}_{\alpha, N}^{N} : \nWN \dashrightarrow \nWN $ by the formulas
\begin{align*}
\hat{\Lambda}_{\alpha, N}^{N+1}(\xx, d\yy)=\mathbf{1}_{\WNNp (\xx)}(\yy)\prod_{k=1}^{N} \hat{m}_\alpha(y_k)  d\yy, \, \hat{\Lambda}_{\alpha, N}^{N}(\xx, d\yy)=\mathbf{1}_{\nWNN (\xx)}(\yy) \prod_{k=1}^{N} \hat{m}_\alpha (y_k) d\yy
.\end{align*}
Introduce sub-Markov Karlin-MacGregor semigroups by the formulas
\begin{align*}
\mathcal{P}_{\alpha,t }^{N}(\xx,d\yy)=\det_{i,j=1}^{N} [p_{\alpha, t}(x_i,y_j)] d\yy, \qquad \hat{\mathcal{P}}_{\alpha,t}^{N}(\xx,d\yy)=\det_{i,j=1}^{N} [\hat{p}_{\alpha, t}(x_i,y_j)] d\yy
.\end{align*}
The semigroups $\mathcal{P}_{\alpha,t }^{N}, \hat{\mathcal{P}}_{\alpha,t}^{N}$  are not necessarily Markov.
Indeed, $\mathcal{P}_{\alpha,t }^{N}$  is associated with $N$-particles $L_\alpha $ diffusion killed when their paths intersect, and similarly $\hat{\mathcal{P}}_{\alpha,t}^{N}$ for $\hat{L}_\alpha $ diffusion. 

\begin{lem}\cite[(13.26)]{AOW19}
Let $\alpha \in\R$.
Then, the following two equations hold for any $t>0$:
\begin{align} \label{:12a}
\mathcal{P}_{\alpha,t}^{N} \hat{\Lambda}_{\alpha, N}^{N}&=\hat{\Lambda}_{\alpha, N}^{N} \hat{\mathcal{P}}_{\alpha,t}^{N}
\\ \label{:12b}
\mathcal{P}_{\alpha,t}^{N+1} \hat{\Lambda}_{\alpha,N}^{N+1}&=\hat{\Lambda}_{\alpha, N}^{N+1} \hat{\mathcal{P}}_{\alpha,t}^{N}
\end{align}
\end{lem}

We first show the shifted intertwining relation \eqref{:34a} using the same technique as in \cite[Theorem 5.1]{Ass20}.
\begin{lem}\label{l:34}
Let $\alpha>-1$.
For any $N\in\mathbb{N} $, $f\in C_{\infty} (\nWN )$, and $t \ge 0$, we have
\begin{align*}
T_{\alpha,t}^{N+1} \LNNp  f=\LNNp T_{\alpha +1,t}^{N}f
.\end{align*}
\begin{proof}
Since this conclusion is clear for $t=0$, we assume $t>0$.
Furthermore, we first consider the case $\xx\in \intnWN$.
Recall $\lambda ^{N}=- N(N-1)/2$.
Multiplying both sides of \eqref{:12b} by 
\begin{align*}
 e^{-\lambda _{}^{N+1} t }  \frac{N! \Delta_N (\yy )}{\Delta _{N+1} (\xx )} \prod_{k=1}^N (\hat{m}_{\alpha}^{} (y_k) )^{-1}
,\end{align*}
we see that the left hand side becomes
\begin{align*}&
\Big( \int_{ }  \frac{N! \Delta_N( \yy )}{\Delta _{N+1} (\zz )} \mathbf{p}_{\alpha,t}^{N+1}(\xx, \zz)  \mathbf{1}_{ \WNNp (\zz ) }(\yy)  d\zz \Big) d\yy 
=(T_{\alpha,t}^{N+1} \LNNp )(\xx,d\yy )
.\end{align*}
On the other hand, the right hand side becomes 
\begin{align*}&
\Big( \int_{ }  \mathbf{1}_{  \WNNp (\xx ) }(\zz) \det_{i,j=1}^N  [ p _{\alpha+1, t}^{}(z_i,y_j)]e^{-\lambda_{}^N t} \frac{N! \Delta _N(\yy )}{\Delta _{N+1} (\xx )}  d\zz \Big) d\yy 
\\&
=\Big( \int_{ }  \mathbf{1}_{ \WNNp (\xx) }(\zz ) \mathbf{p}_{\alpha+1 ,t}^{N} (\zz , \yy) \frac{N! \Delta _N(\zz )}{\Delta _{N+1} (\xx )}   d\zz\Big)d\yy 
\\&
=(\LNNp T_{\alpha+1, t}^{N}  )(\xx,d\yy)
.\end{align*}
Here, we use the fact $\lambda ^{N+1}=N (-1) + \lambda ^{N}$ and \eqref{:21c}.
Combining these, we obtain 
\begin{align*}
T_{\alpha,t}^{N+1} \LNNp (\xx, d\yy) =\LNNp T_{\alpha+1,t}^{N} (\xx,d\yy) 
,\end{align*}
which implies $T_{\alpha,t}^{N+1} \LNNp  f(\xx) =\LNNp T_{\alpha +1,t}^{N}f (\xx ) $ for $\xx\in \intnWN$.
We can extend this for $\xx\in \nWN$ because of the Feller property \eqref{:24o} and \lref{l:33}.
Thus we finish the proof.
\end{proof}
\end{lem}

The Laguerre processes $T_{\alpha+1,t}^{N} $ and $T_{\alpha,t}^{N}$ satisfy another shifted intertwining relation through $\LaNN $ as follows.
\begin{lem}\label{l:44}
Let $\alpha >-1$.
For any $N\in\mathbb{N} $, $f\in C_{\infty} (\nWN )$, and $t \ge 0$, we have
\begin{align*}
T_{\alpha+1,t}^{N} \LaNN f=\LaNN T_{\alpha,t}^{N}f
.\end{align*}
\begin{proof}
We can suppose $t>0$, and we consider the case $\xx\in \intnWN$ first.
Let $\mathbf{1}_{N,\xx}(\yy) $ denote the characteristic function $ \mathbf{1}_{ \nWNN (\xx ) }(\yy)$.
Then, by the definition of $\LaNN $ in \eqref{:32c}, the equality
\begin{align}\label{:36a}
T_{\alpha+1,t}^{N} \LaNN (\xx,d\yy)=\LaNN T_{\alpha,t}^{N}(\xx,d\yy)
\end{align}	
is equivalent to the relation
\begin{align}\label{:32b}
\int_{ } d\zz  \mathbf{1}_{ N, \zz }(\yy) \det_{i,j=1}^N[p_{\alpha+1,t}^{}(x_i, z_j)]  \prod_{k=1}^{N} \frac{x_k^{\alpha+1} }{z_k^{\alpha+1}}  =\!\!\int\!\! d\zz \mathbf{1}_{ N, \xx  }(\zz) \det _{i,j=1}^N[p_{\alpha,t}^{}(z_i, y_j)] \! \prod_{k=1}^{N} \frac{ z_k ^{\alpha}}{y_k^{\alpha} }
.\end{align}

We obtain \eqref{:32b} by the following calculation:
\begin{align*}&
\int_{ } d\zz \mathbf{1}_{N,\zz}(\yy) \det_{i,j=1}^N[p_{\alpha+1,t}^{}(x_i, z_j)]  \prod_{k=1}^{N} \frac{x_k^{\alpha+1} }{z_k^{\alpha+1}}  
\\&&&
\\&
= e^{(\alpha+1)Nt} \int_{ } d\zz  \mathbf{1}_{N,\zz}(\yy) \det_{i,j=1}^N[p_{-\alpha-1,t}^{}(x_i, z_j)] && \qquad\text{ from \eqref{:32a}}
\\&
=e^{(\alpha+1)Nt}  \int_{ } d\zz \mathbf{1}_{N,\zz}(\yy) \det_{i,j=1}^N [\hat{p}_{-\alpha-1 ,t}^{}(z_i,y_j)] \prod_{k=1}^N \frac{ \hat{m}_{-\alpha-1}   (z_k)}{\hat{m}_{-\alpha-1}   (y_k)}  &&  \qquad\text{ from \eqref{:12a}}
\\&
=e^{\alpha Nt} \int_{ } d\zz \mathbf{1}_{N,\zz}(\yy) \det_{i,j=1}^{N}  [p_{-\alpha ,t}^{}(z_i,y_j)]
&& \qquad\text{ from \eqref{:21c}}
\\&
= \int_{ } d\zz \mathbf{1}_{N,\zz}(\yy) \det_{i,j=1}^N [p_{\alpha ,t}^{}(z_i,y_j)]  \prod_{i=1}^{N} \frac{z_k^{\alpha}  }{y_k^\alpha}&& \qquad\text{ from \eqref{:32a}}
.\end{align*}
Thus, we have proved \eqref{:36a} for $\xx\in \intnWN$.
From the Feller properties established in \pref{p:24} and \lref{l:33}, we can extend \eqref{:36a} to all $\xx\in \nWN$, which completes the proof.
\end{proof}
\end{lem}

\begin{remark}
It would be interesting to establish a direct conceptual connection between the intertwining property and the Gibbs property studied by Benigni, Hu and Wu (``Determinantal structures for Bessel fields'', arXiv:2109.09292).
In particular, Lemma 3.3 in their work is a result for the non-colliding square Bessel processes similar to \lref{l:44}.
\end{remark}

\medskip
\noindent
{\em Proof of \tref{t:21}  } \quad 
From \lref{l:34} and \lref{l:44}, we have
\begin{align*}
T_{\alpha, t}^{N+1} \LNNp \LaNN = \LNNp T_{\alpha+1, t}^{N} \LaNN= \LNNp \LaNN T_{\alpha, t}^{N} 
.\end{align*}
Therefore, the equality $\LNNp \LaNN =\LaNNp $ established in \pref{p:23} concludes the statement of this theorem.
\qed 

\begin{cor}
For any $\alpha>-1$, we have \eqref{:26b}.
\begin{proof}
For any $N\in \mathbb{N} $, it holds that $m_\alpha^{N} $ is the unique invariant measure with respect to $T_{\alpha, t}^{N}$ \cite[Remark 1.11]{Ass19b}.
In other words, if a probability measure $m $ on $\nWN $ satisfies $m = m T_{\alpha,t}^{N} $, then $m $ is identical to $m_\alpha^N$.
Now, using the invariance of $m_\alpha^{N+1}$ with respect to $T_{\alpha,t}^{N+1}$ and \tref{t:21}, we have
\begin{align*}
m_\alpha^{N+1} \LaNNp =m_\alpha^{N+1} T_{\alpha,t}^{N+1} \LaNNp = m_\alpha^{N+1}  \LaNNp T_{\alpha,t}^{N}
.\end{align*}
Thus, we obtain $m_\alpha^{N+1} \LaNNp =m_\alpha^{N} $.
\end{proof}
\end{cor}

\section{Proof of \tref{t:12}}\label{s:4}
In this section, $\alpha $ is supposed to be a non-negative integer.
We begin by explaining an interpretation,  established in \cite{KKS16}, of the measure $\LaNN $ defined in \eqref{:32c} in the context of random matrix theory.
Recall that $P_{\mathfrak{rad}}^n [X]  $ denotes the distribution of the radial part of a random matrix $X \in M_{m,n}(\C) $ (cf. \cite{Buf15, Buf16a, Buf16b}).
\begin{lem} \label{l:15}
Assume that $\alpha$ is a non-negative integer.
Let $V_{N+\alpha +1} \in \mathbb{U}(N+\alpha+1)$ be a Haar distributed unitary matrix.
Then, for any $\zz= (z_1 ,\ldots, z_N) \in \nWN $, we have
\begin{align}\label{:33b}
\LaNN (\zz ,\cdot)= P_{\mathfrak{rad}}^{N} [ \pi_{N+\alpha, N}^{N+\alpha+1}( V_{N+\alpha+1} ) \mathrm{diag}(\sqrt{z_1},\ldots, \sqrt{z_{N}}) ] 
.\end{align}

\begin{proof}
It is sufficient to show \eqref{:33b} for $\zz\in\intnWN$.
Applying \cite[Theorem 2.1]{KKS16} for the setting $(m,l, n, \nu)=(N+\alpha+1, N, N,\alpha)$, we have the following result: for $\zz\in \intnWN $, there exists a probability density of 
\begin{align*}
P_{\mathfrak{rad}}^{N} (\pi_{N+\alpha, N}^{N+\alpha+1}( V_{N+\alpha+1} )  \mathrm{diag}(\sqrt{z_1},\ldots, \sqrt{z_{N}})  )
\end{align*}
on $[0,\infty)^N$, and it is proportional to
\begin{align*}
\bigg( \prod_{k=1}^{N} \frac{y_k^{\alpha} }{z_k^{\alpha+1}} \bigg) \frac{\Delta_{N}(\yy)}{\Delta_{N}(\zz)} \det _{i,j=1}^{N} [(z_i-y_j)_{+}^0]
.\end{align*}
Thus, noting that $\det_{i,j=1}^{N}[(z_i-y_j)_{+}^0]=\det_{i,j=1}^{N}[\mathbf{1}_{z_i-y_j\ge 0}]=\mathbf{1}_{X_\ge^{N,N}(\zz )}(\yy)$ for $\zz ,\yy \in \nWN$, we obtain \eqref{:33b} for $\zz\in\intnWN$. 
\end{proof}
\end{lem}

\medskip
\noindent
{\em Proof of \tref{t:12} } \quad
Setting $X_{N+\alpha+1,N} =\pi_{N+\alpha+1, N}^{N+\alpha+1, N+1} (X_{N+\alpha+1,N+1} )$, we have $$\pi_{N}^{N+1}( X_{N+\alpha+1,N+1}^{*} X_{N+\alpha+1, N+1} )=X_{N+\alpha+1,N}^{*} X_{N+\alpha+1,N}.$$
Furthermore, because the Hermitian matrix $ X_{N+\alpha+1,N+1}^{*} X_{N+\alpha+1, N+1} $ is $\mathbb{U}(N+1)$-invariant by conjugation, \eqref{:31c} yields
\begin{align}\label{:12c}
P_{\mathfrak{rad}}^{N+1} [ X_{N+\alpha+1,N+1} ]  \LNNp =  P_{\mathfrak{rad}}^{N} [ X_{N+\alpha+1,N} ] 
.\end{align}

For a random variable $(z_1,\ldots, z_N)$ distributed as $P_{\mathfrak{rad}} ^N[ X_{N+\alpha+1,N}] $, we set 
\begin{align*}
D_{N+\alpha+1,N}=
\begin{bmatrix}
D_{N}  \\
\mathbf{0}_{(\alpha+1 )\times N}
\end{bmatrix}
, D_{N}=\mathrm{diag}(\sqrt{z_1},\ldots,\sqrt{z_N})
.\end{align*}
Let $U_{N}\in \mathbb{U} (N)$ and $V_{N+\alpha+1}\in \mathbb{U} (N+\alpha+1)$ be Haar distributed random matrices such that $D_{N+\alpha+1,N}, U_{N}$, and $ V_{N+\alpha+1}$ are independent.
Then, by a similar reason as in \cite[Lemma 2.4]{Def10}, we have
\begin{align}\label{:12e}
X_{N+\alpha+1, N} \stackrel{\mathrm{law}}{=} V_{N+\alpha+1} D_{N+\alpha+1,N} U_N
.\end{align}
Actually, we can write $X_{N+\alpha+1, N} =V D_{N+\alpha+1, N} U$ with $U \in \mathbb{U}(N)$ and $V \in \mathbb{U}(N+\alpha +1)$.
Because $X_{N+\alpha+1,N}$ is $\mathbb{U} (N+\alpha+1 ) \times \mathbb{U} (N )$-invariant, $X_{N+\alpha+1, N}$ has the same distribution as $(\tilde{V} V) D_{N+\alpha+1, N} (U\tilde{U})$, where $\tilde{U} \in \mathbb{U}(N)$ and $\tilde{V} \in \mathbb{U}(N+\alpha +1)$ are Haar distributed random matrices independent ot $X_{N+\alpha+1,N}$.
The Haar measure being invariant by multiplication, we obtain \eqref{:12e}.

A straightforward computation with \eqref{:12e} yields
\begin{align*}
X_{N+\alpha, N} =\pi_{N+\alpha, N}^{N+\alpha+1, N} (X_{N+\alpha+1, N}) \stackrel{\mathrm{law}}{=} \pi_{N+\alpha,N}^{N+\alpha+1} (V_{N+\alpha+1} ) D_{N} U_N
.\end{align*}
Therefore, using this with the expression of $\Lambda_{\alpha,N}^{N}$ as \eqref{:33b},  we have
\begin{align}\label{:12d}
P_{\mathfrak{rad}}^{N} [X_{N+\alpha, N}]  = P_{\mathfrak{rad}}^{N } [\pi_{N+\alpha,N}^{N+\alpha+1} (V_{N+\alpha+1} ) D_{N} U_{N}]= P_{\mathfrak{rad}}^{N} [ X_{N+\alpha+1,N} ] \Lambda_{\alpha, N}^{N}
.\end{align}
Collecting \eqref{:12c} and \eqref{:12d} with \pref{p:23}, we complete the proof of \tref{t:12}.
\qed

\begin{remark}
If $X_{N+\alpha+1,N} $ is the Wishart matrix, then the relation \eqref{:12d} follows from \cite[Proposition 4.8]{Def10}.
Similar results for generalised Wishart matrices are obtained in \cite{DiW09}.
See also \cite{AMW13} and ``Matrix models for multilevel Heckman-Opdam and multivariate Bessel measures'' by Sun  (arXiv:1609.09096).
\end{remark}

\section*{Acknowledgement}
K.Y. is supported by JSPS KAKENHI Grant Number JP21K13812.


\end{document}